# NEIGHBORING CLUSTERS IN BERNOULLI PERCOLATION[1]

By Ádám Timár

*Indiana University*

We consider Bernoulli percolation on a locally finite quasi-transitive unimodular graph and prove that two infinite clusters cannot have infinitely many pairs of vertices at distance 1 from one another or, in other words, that such graphs exhibit "cluster repulsion." This partially answers a question of Häggström, Peres and Schonmann.

**1. Introduction.** We shall consider Bernoulli($p$) bond percolation on some quasi-transitive locally finite unimodular graph $G$.

Given different infinite clusters $C$ and $C'$, let $\tau(C, C')$ be the set of vertices in $C$ that have distance 1 from $C'$. Call $\tau(C, C')$ the set of *touching points* for $C'$ in $C$. If two infinite clusters touch each other in infinitely many vertices, we say that they are *infinitely touching*; otherwise, they are *finitely touching* (this includes the case where they do not touch at all).

Häggström, Peres and Schonmann asked if there exists some quasi-transitive $G$ and a $p$ such that there are infinitely touching clusters at Bernoulli($p$) percolation. In [6], these authors say that a graph exhibits *cluster repulsion* at level $p$ if two clusters always touch finitely for Bernoulli($p$) bond percolation. They mention that for any $G$, there can be at most countably many such values of $p$ where there is no cluster repulsion. We shall prove that a quasi-transitive unimodular graph always exhibits cluster repulsion. Our proof can be adapted to site percolation without any difficulty. The case of nonunimodular graphs is still open. To avoid meaningless cases, we assume that $p$ is such that the percolation has infinitely many infinite clusters.

THEOREM 1.1. *Let $G$ be a quasi-transitive unimodular graph and consider Bernoulli($p$) edge percolation on it. Then any two infinite clusters touch each other in only finitely many vertices almost surely.*

Received October 2005; revised May 2006.

[1]Supported in part by NSF Grant DMS-02-31224 and Hungarian National Foundation for Scientific Research Grant TO34475.

*AMS 2000 subject classifications.* Primary 60K35, 82B43; secondary 60B99.

*Key words and phrases.* Cluster repulsion, percolation, nonamenable, touching clusters.

Geometric properties of a percolation and how these are related to certain properties of the underlying graph have been intensively studied. For an overview and references, see [7]. Some of the fundamental results connected to our subject are now outlined. It is conjectured that a transitive graph can have infinitely many infinite clusters at Bernoulli($p$) percolation for some $p$ if and only if the graph is nonamenable. It is known that infinitely many infinite clusters imply nonamenability; the converse was shown for *some* Cayley graph of an arbitrary group. It is well known that if there are infinitely many infinite clusters, then each of them has infinitely many ends. Moreover, still assuming infinitely many infinite clusters, there are infinitely many points having the property that deletion from its cluster results in at least three infinite components. Such points are called *encounter points* and were introduced in [4]. Encounter points with particular properties will play an important role in our proof.

We note that even for amenable graphs, there are group-invariant random subgraphs where any two infinite clusters touch infinitely. An example is the uniform spanning forest $\mathbb{Z}^d$ for $4 < d \leq 8$; see [2]. However, our arguments can be applied to any insertion and deletion tolerant percolation. (The definitions are given later in this section; the transfer from Bernoulli to these more general percolations is explained in Remark 2.6.) On the other hand, if we do not assume quasi-transitivity, then there is an example for a graph with infinitely touching clusters at some Bernoulli percolation, as claimed in [6].

An essential tool in the study of group-invariant percolations is the so-called Mass Transport Principle (MTP). A corresponding graph property is unimodularity: we say that a quasi-transitive graph is *unimodular* if there is a $K > 0$ such that for any two $x, y \in V(G)$, we have $|S_x y|/|S_y x| \leq K$. Here, $S_x$ is the stabilizer of $x$ in the group of automorphisms of $G$ and $S_x y$ denotes the orbit of $y$ by the stabilizer of $x$. Every Cayley graph is transitive and unimodular.

The MTP was first used in percolation theory in [5] and was developed more generally in [3]. We state here a simple corollary of the principle.

PROPOSITION 1.2. *Let $G$ be a connected locally finite graph whose group of automorphisms is unimodular and acts quasi-transitively. If $T(x, y)$ is a nonnegative invariant function [i.e., $T(x, y) = T(gx, gy)$ for every automorphism $g$], then $\sum_y T(x, y)$ is finite for every $x \in V(G)$ whenever $\sum_y T(y, z)$ is finite for every $z \in V(G)$.*

Usually, we shall define some group-invariant random function in two variables, determined by the configuration of the percolation. A value $t(x, y) = c$ of such a function will be defined by saying "let $x$ send mass $c$ to $y$." Then $T(x, y)$ will be the expectation of $t(x, y)$ over all the random configurations;

we shall refer to it as the *expected mass* sent from $x$ to $y$. Sometimes, we do not specify $t(x,y)$ for every pair $(x,y) \in V(G)^2$—in these cases, every $t(x,y)$ not defined is automatically 0.

Hereafter, we shall always assume that there is a Bernoulli($p$) edge percolation on an underlying locally finite unimodular quasi-transitive graph $G$, without always mentioning these assumptions. Let $d$ be the maximal degree in $G$.

A process is said to be *insertion tolerant* if for any edge $e$ and any event $A$, $\mathbf{P}[A] > 0$ implies $\mathbf{P}[\{\kappa \cup e : \kappa \in A\}] > 0$. When creating the event $\{\kappa \cup e : \kappa \in A\}$ from $A$, we say that $e$ was *inserted*. Define *deletion tolerance* analogously. Note that Bernoulli percolation is insertion and deletion tolerant.

We can typically use insertion tolerance in the following setting. Consider some property that a finite set of edges may have in a configuration. For a finite subset $X$ of edges, let $A(X)$ be the event that $X$ has the property and let $A$ be the event that some finite $X$ has the property. Since $A = \bigcup_X A(X)$, and this is a countable union, we have that $\mathbf{P}[A] > 0$ implies $\mathbf{P}[A(X)] > 0$ for some $X$. Hence, when we show that some configuration contains an $X$ with the property and say "insert $X$," that is a shorthand way of saying that we take an $X$ that satisfies the property with positive probability and insert $X$ on this event. Similar terminology is applied for deletion.

Say that two vertices (sets of vertices) *$k$-touch* each other if their distance in $G$ is at most $k$. Two clusters *$k$-touch in infinitely many points* if one of them contains infinitely many points that $k$-touch the other.

Given clusters $C$ and $C'$, denote by $\tau(C, C')$ the set of vertices in $C$ that touch $C'$. The *number of ends in $\tau(C, C')$* is defined to be the supremum of the number of infinite components of $C \setminus F$ that contain infinitely many points of $\tau(C, C')$ as $F$ ranges through finite subsets of $C$.

In the course of the proof, we shall refer to mass transports according to the following scheme. Suppose that there is some automorphism-equivariant function $f$ that assigns a finite nonempty set $f(C, C')$ of vertices to certain pairs of infinitely touching clusters $C$ and $C'$. For each vertex $x$, consider the cluster $C$ of $x$ and if $x$ touches $C'$, then let $x$ send mass $1/|f(C, C')|$ to each element of $f(C, C')$. The expected mass sent out is at most $d$, while the expected mass received is infinite if some $\tau(C, C')$ is infinite, a contradiction to Proposition 1.2. So, there cannot exist such a function $f$ if $\tau(C, C')$ is infinite. Using this scheme, we are able to use standard mass transport arguments that were developed for clusters (and which do not generally work for pairs of clusters). A general corollary of the Mass Transport Principle (MTP) for unimodular graphs is that one cannot assign a finite set of vertices to each infinite cluster in some equivariant way. Of course, this is not generally true for pairs of infinite clusters (consider, e.g., the set of touching points between pairs of finitely touching clusters). However, we can extend it to pairs of infinitely touching clusters by means of the above argument, that is, one cannot assign a finite set of vertices to each pair of infinitely touching clusters in an equivariant way.

Before proceeding to the proof, we give an overview of it. Suppose that there exist infinite $\tau(C, C')$'s. The first step is to show that some $\tau(C, C')$ has infinitely many ends, as follows. One can use insertion and deletion tolerance to show that there are $\tau(C, C')$'s with $\geq 2$ ends. The argument is based on the fact that we may assume $G$ is 2-connected and can then choose two disjoint paths between the elements of two pairs of infinitely touching clusters. Surprisingly, the proof that there are $\tau(C, C')$'s with $\geq 3$ ends is not that straightforward. Now, if $\tau(C, C')$ had only two ends, then some of these $\tau(C, C')$'s have infinitely many ("good") cutedges $e$ with the property that $C \setminus e$ has two components with infinitely many elements of $\tau(C, C')$. (There will be some extra technicalities required of a good edge.) The existence of good edges will be provided by our construction of 2-ended $\tau(C, C')$'s. Fix an $o$ and $o'$ and let $E$ be the event that for the cluster $C$ of $o$ and the cluster $C'$ of $o'$, $\tau(C, C')$ has exactly two ends. Let $e_1, e_2, \ldots$ and $e_{-1}, e_{-2}, \ldots$ be the good cutedges, in the directions of the ends of $\tau(C, C')$ (resp. in the order of their distances from $o$ in $C$). We define an event $E_i$ from $E$ by closing the edges $e_i$ and $e_{-i}$ in each configuration. It will turn out that $\mathbf{P}[E_i] \geq c\mathbf{P}[E]$ for some constant $c$. On the other hand, any configuration is contained in at most a bounded number of the $E_i$'s, hence $\inf_i \mathbf{P}[E_i] = 0$. This implies $\mathbf{P}[E] = 0$ and we conclude that $\tau(C, C')$ has $\geq 3$ ends with positive probability. Then $\tau(C, C')$ has infinitely many ends, by the MTP used for infinitely touching clusters. In particular, $\tau(C, C')$ will have exponential growth in $C$ with positive probability (meaning that its elements intersect the ball of radius $r$ in $\beta^r$ elements for some $\beta > 1$ and $r$ sufficiently large). Finally, we shall define a mass transport. For each vertex $x$ and each of its neighbors $y$, choose a random minimal path between $x$ and $y$ (in their cluster if they are in the same cluster) and let $x$ send mass $1, 1/2^2, 1/3^2, \ldots$ to the consecutive vertices on this path, starting from $x$. The expected mass sent out is finite, but the expected mass received is infinite, as shown by the endpoints of an edge inserted between clusters $C$ and $C'$ with $\tau(C, C')$ having exponential growth.

**2. No infinite touchings.** First, we are going to show that there are $\tau(C, C')$'s with infinitely many ends. We need two simple graph-theoretic lemmas.

LEMMA 2.1. *Let $G$ be some 2-connected quasi-transitive graph. Then any finite subgraph $H$ of $G$ is contained in some 2-connected finite subgraph $H'$ of $G$.*

PROOF. If $H$ does not have any cutvertices, then $H' := H$ immediately gives the claimed assertion. Otherwise, for each cutvertex $x \in H$ and each pair $A, B$ of components of $H \setminus x$, choose an arbitrary path $P(A, B)$ in $G$ that joins $A$ and $B$ and does not contain $x$. Such a path exists by the 2-connectivity of $G$. Let $H'$ be the union of $H$ and the set of paths $P(A, B)$ over all choices of $x, A$ and $B$. If $H'$ had a cutvertex $x$, then only one component of $G \setminus x$ could contain any vertex from $H$ [otherwise, there would have been a path not containing $x$ between vertices of two such components: either a path in $H$ or a $P(A, B)$]. So, all but at most one component of $G \setminus x$ consists only of vertices in $P(A, B)$'s. Let $D$ be one such component. The endpoints of every $P(A, B)$ are in $H$, so at least one of them differs from $x$ and is in $D \cap H$, giving a contradiction. □

A 2-*connected component* (block) of a graph $G$ is a maximal 2-connected subgraph of $G$.

LEMMA 2.2. *If $G$ is some infinite graph whose automorphism group is unimodular and acts quasi-transitively, then for any 2-connected component $C$ of $G$, the stabilizer of $C$ in $\mathrm{Aut}(G)$ acts quasi-transitively on $C$.*

PROOF. The claim is trivial if the 2-connected components are finite. The infinite 2-connected components of $G$ are the infinite components of the graph $G'$, constructed as follows. Let the vertex set $V(G')$ consist of all pairs $(x, C)$, where $x \in V(G)$ and $C$ is the 2-connected component of $G$ that contains $x$. Two vertices $(x, C)$ and $(x', C')$ are defined to be *adjacent* in $G'$ if $C = C'$ and $\{x, x'\}$ is an edge of $G$.

Another way of constructing $G'$ from $G$ is as follows. Suppose that $G \setminus v$ has components $C_1(v), C_2(v), \ldots, C_k(v)$ for vertex $v \in V(G)$. Replace every vertex $v$ by vertices $v_1, \ldots, v_k$ and connect $v_i$ to $u_j$ if $v$ and $u$ are adjacent in $G$, $v \in C_i(u)$ and $u \in C_j(v)$. Now, if $G$ had cutvertices, then we obtain a graph $G'$ with infinitely many components. If two edges are contained in a cycle of $G$ (which holds iff they are in the same 2-connected component), then the corresponding edges in $G'$ are also contained in a cycle and hence they are in the same component of $G'$. It is clear that $\mathrm{Aut}(G)$ acts quasi-transitively on $G'$: the transitivity class of $v_i$ is determined by the transitivity class of $v$ and by $i$. Hence, the stabilizer in $\mathrm{Aut}(G)$ of any component of $G'$ is quasi-transitive on the component. □

By Lemma 2.2, we may assume in what follows that $G$ is 2-connected since all elements of a $\tau(C, C')$ are in the same 2-connected component.

LEMMA 2.3. *Suppose that there exist infinitely touching clusters. Then almost always there are clusters $C$ and $C'$ such that $\tau(C, C')$ has infinitely many ends.*

PROOF. It is enough to show that some $\tau(C, C')$ has at least three ends. It then follows by deletion tolerance that for some $\tau(C, C')$, there are vertices $x \in C$ such that $C \setminus x$ has at least three infinite components with infinitely many elements of $\tau(C, C')$ in each. Then the standard MTP argument can be applied to show the existence of infinitely many such vertices and hence the existence of infinitely many ends [since the number of $\tau(C, C')$'s that a vertex can be contained in is bounded].

We first prove that there exist different components $C_1, D_1, C_2, D_2$ with the following properties: $C_i$ infinitely touches $D_i$ and finitely touches $C_{3-i}$ and $D_{3-i}$, similarly for $C$ and $D$ interchanged. To show this, define a random graph $H$ whose vertices are the infinite clusters of $G$ and put an edge between two if they infinitely touch. We need there to be four vertices in $H$ that induce two disjoint edges. Then the clusters corresponding to the endpoints of these edges will supply $C_1, C_2, D_1, D_2$. Suppose that, on the contrary, for any pair of disjoint edges in $H$, there is some other edge joining the four endpoints. This can be used to show that $H$ contains vertices $x$ and $y$ such that there are three disjoint paths between $x$ and $y$. Namely, choose an arbitrary (self-avoiding) path of length 8 in $H$ with consecutive edges $e_1, \ldots, e_8$. (Note that each vertex in $H$ has the same degree, by indistinguishability of infinite clusters, as proved in [8]. Hence, a path of length 8 exists in $H$, or otherwise $H$ would only contain finite components. These finite components are isomorphic graphs, again by indistinguishability, and they are not empty, by our assumption. Choosing two edges from different components would provide us with four vertices that induce two disjoint edges.) Now, we assumed that the endpoints of $e_1$ and $e_8$ induce some other edge $f_1$ in $H$ and that the endpoints of $e_3$ and $e_6$ also induce some other edge $f_2$ in $H$. Then the endpoints $x$ and $y$ of $f_2$ can be connected by three disjoint paths in $H$: a path containing $e_4$ and $e_5$ (possibly among others) a path containing $e_2, f_1$ and $e_7$ (possibly among others) and, finally, the path consisting of $f_2$. Let the vertices on the first of these paths be $x, q_1, \ldots, q_i, y$. The corresponding infinite clusters $X, Q_1, \ldots, Q_i, Y$ in $G$ can be used to find a path from $X$ to $Q_i$ that only intersects $X, Q_1, \ldots, Q_i$. Insert it. Do this with the other two paths in $H$ too. We get a cluster $K$ such that $\tau(K, Y)$ has at least three ends, so the claimed assertion is proved. We conclude that we may assume that there exist clusters $C_1, C_2, D_1, D_2$ with the above properties with probability 1.

By deletion tolerance, we may also assume that for any $X, Y \in \{C_1, C_2, D_1, D_2\}$, if $X$ and $Y$ finitely touch, then $\text{dist}(X, Y) \geq 2$. (If the distance is 1, delete the finitely many vertices in $X$ that touch $Y$. This may break $X$ into finitely many pieces, but one of them will still infinitely touch the cluster in $\{C_1, C_2, D_1, D_2\}$ that $X$ infinitely touched and will finitely touch the other two.)

In what follows, we are going to delete finitely many edges and insert two disjoint paths between the elements of two disjoint pairs of $\{C_1, C_2, D_1, D_2\}$. Finally, we get two clusters $C, C'$ with a 2-ended $\tau(C, C')$.

Let $\Gamma$ be a finite 2-connected subgraph of $G$ that intersects all of $C_1, C_2$, $D_1, D_2$. Such a choice exists by Lemma 2.1. (More precisely, $\Gamma$ intersects a 4-tuple of such clusters with positive probability. However, we will continue to use the language introduced in Section 1, concerning the equivalence of finding a finite subgraph with a certain property for a particular configuration and finding a finite subgraph that satisfies this property with positive probability.) Now, let $C'_1$, $C'_2$, $D'_1$, $D'_2$ be clusters of $G \setminus \Gamma$ contained in $C_1$, $C_2$, $D_1$, $D_2$, respectively, such that $C'_i$ infinitely touches $D'_i$ (and finitely touches the other two), similarly for $D'_i$. Let $c_i \in \Gamma$ (resp. $d_i \in \Gamma$) be a vertex that is adjacent to $C'_i$ (resp. $D'_i$) in $G$.

Since $\Gamma$ is 2-connected, there exist vertex-disjoint paths $Q$ and $Q'$ in $\Gamma$ that connect $c_1$ and $d_1$, respectively to the set $\{c_2, d_2\}$. (It is well known that such paths exist in a 2-connected graph. Add extra vertices $v$ and $w$ to $\Gamma$, connect $v$ to $c_1$ and $d_1$ and connect $w$ to $c_2$ and $d_2$. The resulting graph is still 2-connected, hence there are two inner-disjoint paths between $v$ and $w$ by Menger's theorem and these supply the two paths in $\Gamma$ that we need.) If we deleted the edges of $\Gamma$ and then inserted the edges of $Q$ and $Q'$, we would obtain two clusters $C, C'$ with a 2-ended $\tau(C, C')$. However, we shall need $C$ to have cutedges with some extra property and in order to have such cutedges, we may need to modify $Q$ a little. Namely, take a path $P$ as follows. Suppose that $Q$ connects $c_1$ to $c_2$. Let $x$ be the first vertex on $Q$ (counting from $c_1$) at distance 1 from $C'_2$. Let $e$ be an arbitrary edge that connects $x$ to $C'_2$. Finally, denote by $P$ the subpath of $Q$ from $c_1$ to $x$, followed by $e$. Insert $P$ and $Q'$. (When $Q$ connects $c_1$ to $d_2$, do the same with $c_2$ and $C_2$ replaced by $d_2$ and $D_2$, resp., in the previous description.)

On the resulting event of positive probability, $c_1$ is in an infinite component $C$ and $d_1$ is in a component $C'$ such that $\tau(C, C')$ has at least two ends. Furthermore, $C$ contains an edge $e$ with the following properties:

(1) the deletion of $e$ from $C$ results in two infinite connected components $X_1, X_2$ that infinitely touch $C'$;
(2) the only vertex in $X_1$ at distance 1 from $X_2$ is an endpoint $x(e)$ of $e$.

The edge $e$ of the previous paragraph has these properties: property (2) follows from the fact that the two clusters connected by $Q$ had distance $\geq 2$ and from the way in which we defined $x$ there. We will call edges satisfying properties (1) and (2) *good* in $C$ with respect to $C'$. Sometimes, we simply say "good edges"—in these cases, the $C'$ is clear from the context or is fixed. Since we have at least one good edge, the MTP ensures that $C$ contains infinitely many good edges with respect to $C'$. [Otherwise, let every vertex $x$ for each $x \in \tau(C, C')$ send mass $1/2k$ to the endpoints of the good edges, where $k$ is the number of edges in $C$ that are good with respect to $C'$. The expected mass sent out is $\leq d$, while the expected mass received is infinite.] Note that if a good edge $e$ is deleted from cluster $C$ and the resulting components are $X_1$ and $X_2$, as above, then there are at most $d$ edges in $G$ whose insertion would connect $X_1$ and $X_2$. This is so because such an edge has to be incident to the endpoint $x(e)$ of $e$, as in property (2).

We have seen that there are vertices $o$ and $o'$ such that, with positive probability, the cluster $C$ of $o$ infinitely touches the cluster $C'$ of $o'$ and such that $C$ has infinitely many good edges with respect to $C'$. [Hence, $\tau(C, C')$ has $\geq 2$ ends.] Fix such vertices $o$ and $o'$ and call the described event $E$. Hereafter, $C$ and $C'$ stand for the clusters of $o$ and $o'$, respectively. Suppose that a $\tau(C, C')$ has two ends (otherwise, we are done).

For each $i \in \mathbb{Z}^+$, define the following mapping from $E$ onto an event $E_i$. Consider a configuration $\omega$ in $E$. For any good edge $e$, there are exactly two good edges that can be connected to $e$ in $C$ by paths that do not contain any other good edge [otherwise, $\tau(C, C')$ would have more than two ends]. These two good edges are separated from each other in $C$ by $e$. Thus, for any two disjoint paths in $C$ that contain infinitely many good edges and that start from $o$, the sequence of good edges $e_1, e_2, \ldots$ and $e_{-1}, e_{-2}, \ldots$ along the paths is unique up to the orders of the two paths (i.e., the signs of the indices). Let $\phi_i(\omega)$ denote the configuration that results from $\omega$ if we close $e_i$ and $e_{-i}$. Note that $\phi_i$ is measurable for every $i$. Define the set of resulting configurations to be $E_i$.

Denote the component of vertex $x$ in some configuration $\omega$ by $C_\omega(x)$. We claim that from a configuration $\omega$ in $\bigcup_{i \in \mathbb{Z}^+} E_i$, we can recover $\{x_i, x_{-i}\}$, where $x_i$ (resp. $x_{-i}$) is an endpoint of $e_i$ (resp. $e_{-i}$). This is true because $C_\omega(o)$ contains only finitely many touching points to $C'$ and is adjacent to only two clusters $K$ and $K'$ that infinitely touch $C'$. [If $C_\omega(o)$ were adjacent to some third component $K''$ that infinitely touches $C'$ in $\omega$, then inserting edges between $C_\omega(o)$ and $K$, $K'$, $K''$ gives a component $\hat{C}$ where $\tau(\hat{C}, C')$ has at least three ends, so we are done.] Moreover, since the edges connecting $C_\omega(o)$ to $K$ and $K'$ were good, there is only one vertex of $C_\omega(o)$ that touches $K$ or there is only one vertex of $K$ that touches $C_\omega(o)$, similarly for $K'$. Let two such vertices be $x_i$ and $x_{-i}$ (for $K$ and for $K'$, resp.).

Now, $e_i$ is incident to $x_i$, thus there are at most $d^2$ possible pairs of edges for $e_i, e_{-i}$ and, equivalently, at most $d^2$ configurations $\omega'$ with $\phi_i(\omega') = \omega$ for some $i$. One consequence of this fact is that any configuration is contained in at most $d^2$ of the $E_i$'s and hence

$$\inf_i \mathbf{P}[E_i] = 0.$$

For any $e, e' \in E(G)$, let $F_i(e, e')$ be the set of configurations in $E$ such that $\{e_i, e_{-i}\} = \{e, e'\}$. Thus, $e$ and $e'$ are open on $F_i(e, e')$. Let $\bar{F}_i(e, e') \subset E_i$ be the set of configurations arising from $F_i(e,e')$ if we close $e$ and $e'$, that is, $\phi_i(F_i(e,e')) = \bar{F}_i(e,e')$. It is clear that $p^{-2}\mathbf{P}[F_i(e,e')] = (1-p)^{-2}\mathbf{P}[\bar{F}_i(e,e')]$. Since $E$ is the disjoint union of the $F_i(e,e')$ [as $e,e' \in E(G)$], we have

$$\mathbf{P}[E] = \sum_{e,e' \in E(G)} \mathbf{P}[F_i(e,e')]$$

$$= \left(\frac{p}{1-p}\right)^2 \sum_{e,e' \in E(G)} \mathbf{P}[\bar{F}_i(e,e')]$$

$$\leq \left(\frac{p}{1-p}\right)^2 d^2 \mathbf{P}[E_i],$$

where the last inequality comes from the fact that the *multi*set $\bigcup_{e,e' \in E(G)} \bar{F}_i(e,e')$ contains each element of $E_i$ at most $d^2$ times. This inequality, together with $\inf_i \mathbf{P}[E_i] = 0$, implies that $\mathbf{P}[E] = 0$.

Thus, there are $\tau(C,C')$'s with infinitely many ends. $\square$

LEMMA 2.4. *With the assumptions of the previous lemma, there exist clusters $C$ and $C'$ and a positive number $\alpha$ such that $\tau(C,C')$ has infinitely many ends and, further, it has infinitely many points $x$ that $\alpha$-touch $C'$ and have the property that at least three infinite components of $C \setminus x$ contain infinitely many elements of $\tau(C,C')$.*

PROOF. By Lemma 2.3, there are finite sets $F$ such that at least three components of $C \setminus F$ contain infinitely many elements of $\tau(C,C')$. By deletion tolerance, there also exists such an $F = \{x\}$. If $\alpha$ is chosen to be sufficiently large, then some $x$ with this property will $\alpha$-touch $C'$. There are infinitely many such $x$'s by the MTP. $\square$

Fix an $\alpha$ as in Lemma 2.4. We shall call points with the property of $x$ in Lemma 2.4 *strong encounter points* with respect to $C'$.

Suppose that $T$ is a tree with edges labeled by some positive integers (lengths). We say that $T$ has *exponential labeled growth* if there exists some constant $c > 1$ such that the number of points at distance $n$ from a fixed vertex is at least $c^n$ for every sufficiently large $n$, where the distance of two adjacent vertices is understood to equal the integer labeling the edge between them. Note that if the labels are bounded, then exponential labeled growth is equivalent to exponential growth of the underlying graph.

LEMMA 2.5. *Suppose that there exist infinite clusters $C$ and $C'$ that touch in infinitely many points. Then there exists a forest $\phi(C,C')$, defined on the strong encounter points of $C$ with respect to $C'$, with the following* *properties. The union $\Phi := \bigcup \phi(C, C')$ is also a forest, with automorphism-invariant law. If the length of each edge in $\phi(C, C')$ is defined to be the distance in $C$ between its endpoints, then some tree in $\Phi$ has exponential labeled growth.*

PROOF. By Lemma 2.4, we have strong encounter points. If $\tau(C, C')$ has strong encounter points, then for any such point $x$ and any infinite component $Y$ of $C \setminus x$ such that $|Y \cap \tau(C, C')|$ is infinite, $Y$ will contain strong encounter points by the MTP [otherwise, let each element of $Y \cap \tau(C, C')$ send mass 1 to $x$]. For an arbitrary $x$, let $(C, C_1), \ldots, (C, C_k)$ be the pairs of clusters such that $x$ is a strong encounter point in $C$ with respect to $C_i$. (So $k$ is at most the size of a ball of radius $\alpha$.) Choose an element $j \in \{1, \ldots, k\}$ uniformly at random and let $f(x) := \tau(C, C_j)$. We do this independently for each $x$ that is a strong encounter point for some pair of clusters. For each pair of clusters $C$ and $C'$, we shall define a forest on $\nu(C, C') := \{x : f(x) = \tau(C, C')\}$. Note that the $\nu(C, C')$ are all disjoint, by definition. If there were to exist a strong encounter point in $C$ with respect to $C'$, then there would be infinitely many, so $\nu(C, C')$ is nonempty. Then by our opening remark, for any $x \in \nu(C, C')$ and any infinite component $Y$ of $C \setminus x$ with $Y \cap \tau(C, C') \neq \varnothing$, $Y$ contains infinitely many strong encounter points with respect to $C'$ and hence $\nu(C, C') \cap Y$ is nonempty a.a.

For any $v \in \nu(C, C')$ and each of the (at least three) components of $C \setminus v$ that contain some element of $\nu(C, C')$, choose uniformly an element of $\nu(C, C')$ that has minimal distance from $v$ in this component. Put a directed edge from $v$ to this vertex. Doing this for every $v \in \nu(C, C')$, we obtain a digraph $\vec{H}(C, C') =: \vec{H}$. Denote by $H$ the graph that results from ignoring the directions of the edges in $\vec{H}$. There may be cycles in $H$, but any two cycles share at most one vertex. [Suppose, on the contrary, that there are two cycles that can share more than one vertex and take their union $J$. Since $J$ is 2-connected, it is in the same component of $C \setminus v$ for any $v \in V(C)$, hence the outdegree of each vertex in the restriction of $\vec{H}$ to $J$ is at most 1. Hence, the average degree in $J$ is $\leq 2$. On the other hand, since $J$ is the union of two intersecting cycles, the average degree is $> 2$, a contradiction.] So, we can delete a uniformly chosen edge from each of the possibly arising, pairwise edge-disjoint cycles in $H$ to obtain a forest $F(C, C')$. Now, let the label of each edge in $F(C, C')$ be the distance of its two endpoints (as vertices in $C$) in $C$.

The family of forests $F(C, C')$ (as $C$ and $C'$ range through all infinite clusters) was constructed in an automorphism-invariant way. The vertex sets of the $F(C, C')$'s are disjoint because the $\nu(C, C')$'s were disjoint. So $\bigcup F(C, C')$ is an invariant forest. Every point in every $F(C, C')$ has degree $\geq 3$. Hence, there is some number $k$ such that the subforest of $\bigcup F(C, C')$ consisting of the edges that have label $\leq k$ also has expected degree $>2$. Denote the restriction of this forest to $F(C, C')$ by $\hat{\phi}(C, C')$. By Theorem 7.2 in [3], this is equivalent to $p_c(\{\bigcup \hat{\phi}(C, C')\}) < 1$. Those $\hat{\phi}(C, C')$'s with $p_c < 1$ have exponential growth by an easy and well-known counting argument. We conclude that $F(C, C')$ also has exponential labeled growth for certain pairs $C, C'$. Choose $\phi(C, C') := F(C, C')$ and their family as $\Phi$. □

PROOF OF THEOREM 1.1. Define the following mass transport. (Here $\alpha$ is defined as in the previous lemma.) If $x$ and $y$ are in the same infinite cluster and their distance in $G$ is at most $\alpha$, take uniformly at random a path $x_1, x_2, \ldots, x_m$ of minimal length in the cluster between them ($x_1 = x$). Let $x$ send mass $1/k^2$ to $x_k$. The expected mass sent out by $x$ is at most the size of its $\alpha$-neighborhood times $\sum 1/k^2$, which is finite.

We are now going to show that the expected mass received is infinite. Let $E$ be the event of positive probability that vertex $o$ is in an infinite cluster $C$, it has a neighbor $o'$ that is in another infinite cluster $C'$ and some tree in $\phi(C, C')$ has exponential growth (with lengths on its edges the same as the distances in $C$). Let $c > 1$ be a number such that $|V(\phi(C, C')) \cap S_n(C)| \geq c^n$ for every sufficiently large $n$, where $S_n(C)$ is the set of points at distance $n$ from $o$ within $C$. Insert the edge between $o$ and $o'$. In the resulting event of positive probability, $o$ receives mass $\geq \gamma \sum_n c^n n^{-2}$ with some constant $\gamma$ because every vertex of $V(\phi(C, C')) \cap S_n(C)$ sends mass $1/n^2$ to $o$. So, the expected mass received is indeed infinite. This contradiction completes the proof. □

REMARK 2.6. The proof remains valid, with small modifications, if we assume only that the percolation is insertion and deletion tolerant. The only property of Bernoulli percolation that we have used is strong insertion and deletion tolerance. This means that there is exists a constant $c > 0$ such that for any event $A$ and edge $e$, the inequality $\mathbf{P}[\{\kappa \cup e : \kappa \in A\}] \geq c\mathbf{P}[A]$ holds, similarly for deletion. However, simple insertion and deletion tolerance is enough (with some extra care), by an argument similar to one suggested by Häggström in [8]. The only part of our proof that does not immediately generalize is in Lemma 2.3, where we defined the mappings from $E$ onto the $E_i$'s. A uniform lower bound there such as $\mathbf{P}[E_i] \geq \frac{(p-1)^2}{d^2 p^2} \mathbf{P}[E]$ no longer necessarily exists. We therefore make the following changes in the proof. As before, $\mathbf{P}$ will denote the measure corresponding to the percolation. For an edge $e \in E(G)$, let $F(e) = F$ be the event that $e$ is good in its cluster with respect to some other cluster. Define the following measure $\mu$ on $F$: for every measurable set $A \subset F$, $\mu(A) := \mathbf{P}[\{\kappa \setminus e : \kappa \in A\}]$. By insertion tolerance, $\mu$ is absolutely continuous with respect to $\mathbf{P}$ because if some $B := \{\kappa \setminus e : \kappa \in A\}$ has positive $\mathbf{P}$-measure ($A \subset F$), then $A$ also has positive $\mathbf{P}$-measure (since a subset of $A$ arises from $B$ by inserting $e$). Consider the Radon–Nikodým derivative $f$ of $\mu$ over $\mathbf{P}$. By deletion tolerance, $f$ is $\mathbf{P}$-almost everywhere positive, hence the set $F_\delta(e) = \{\kappa \in F : f(\kappa) > \delta\}$ has positive measure if $\delta > 0$ is sufficiently small. If a configuration is in $F_\delta(e)$, we say that $e$ is $\delta$-*good* in that configuration. Note that there are finitely many possible values for $\mathbf{P}[F_\delta(e)]$ over $e \in E(G)$ because $G$ is quasi-transitive. Thus, there are also infinitely many $\delta$-good edges with positive probability. By the MTP, there are also infinitely many $\delta$-good edges for some pair of clusters. Let $E$ be the event that the cluster of $o$ has infinitely many $\delta$-good edges with respect to the cluster of $o'$. We can continue with the proof of Lemma 2.3, since now $\mathbf{P}[E] > \frac{\delta^2}{d^2}\mathbf{P}[E_i]$ for every $i$.

QUESTION 2.7. Consider some supercritical Bernoulli percolation on some transitive graph and consider adjacent vertices $x$ and $y$. Conditioned on $x$ and $y$ being in the same cluster, we conjecture that their expected distace within the cluster is finite. The distance within the cluster, also called *chemical distance*, was proved to have an exponential decay in case of supercritical percolation on $Z^d$ [1], hence the conjecture is true in that setting. If our conjecture were true for nonamenable graphs and $p_c < p < p_u$, then the result of the present paper would follow from a much simpler argument similar to the final proof of Theorem 1.1.

**Acknowledgments.** I would like to thank Jan Swart, who first told me about the problem, and Olle Häggström for reading the manuscript. I am grateful to the anonymous referees for pointing out a gap and for other comments. I am indebted to Russell Lyons for all his questions and suggestions on previous drafts.

DEPARTMENT OF MATHEMATICS
INDIANA UNIVERSITY
BLOOMINGTON, INDIANA 47405-5701
USA
E-MAIL: atimar@indiana.edu
URL: http://mypage.iu.edu/~atimar/